\documentclass[preprint,12pt,compress]{elsarticle}

\usepackage{amssymb}
\usepackage{amsmath}
\usepackage[a4paper, body={14.8cm,23cm}]{geometry}

\begin{document}
\journal{ArXiv}
\parindent 15pt
\parskip 4pt


 \newcommand{\eps}{\varepsilon}
 \newcommand{\lam}{\lambda}
 \newcommand{\To}{\rightarrow}
 \newcommand{\as}{{\rm d}P\times {\rm d}t-a.e.}
 \newcommand{\ps}{{\rm d}P-a.s.}
 \newcommand{\jf}{\int_t^T}
 \newcommand{\tim}{\times}

 \newcommand{\F}{\mathcal{F}}
 \newcommand{\E}{\mathbf{E}}
 \newcommand{\N}{\mathbb{N}}
 \newcommand{\s}{\mathcal{S}}
 \newcommand{\M}{{\rm M}}
 \newcommand{\T}{[0,T]}
 \newcommand{\LT}{L^1(\R)}

 \newcommand{\R}{{\mathbb R}}
 \newcommand{\Q}{{\mathbb Q}}
 \newcommand{\RE}{\forall}

\newcommand{\gET}[1]{\underline{{\mathcal {E}}_{t,T}^g}[#1]}
\newcommand {\Lim}{\lim\limits_{n\rightarrow \infty}}
\newcommand {\Limk}{\lim\limits_{k\rightarrow \infty}}
\newcommand {\Limm}{\lim\limits_{m\rightarrow \infty}}
\newcommand {\llim}{\lim\limits_{\overline{n\rightarrow \infty}}}
\newcommand {\slim}{\overline{\lim\limits_{n\rightarrow \infty}}}
\newcommand {\Dis}{\displaystyle}

\begin{frontmatter}
\title {$L^1$ solutions to one-dimensional BSDEs with sublinear growth generators in $z$\tnoteref{fund}}
\tnotetext[fund]{Supported by the National Natural Science Foundation of China (No. 11371362), the China Postdoctoral Science Foundation (No. 2013M530173 and 2014T70386), the Qing Lan Project and the Fundamental Research Funds for the Central Universities (No. 2013RC20)\vspace{0.2cm}.}

\author{ShengJun FAN}
\ead{f$\_$s$\_$j@126.com}

\address{School of Mathematics, China University of Mining and Technology, Xuzhou 221116, PR China\vspace{-0.8cm}}

\begin{abstract}
This paper aims at solving a one-dimensional backward stochastic differential equation (BSDE for short) with only integrable parameters. We first establish the existence of a minimal $L^1$ solution for the BSDE when the generator $g$ is stronger continuous in $(y,z)$ and monotonic in $y$ as well as it has a general growth in $y$ and a sublinear growth in $z$. Particularly, the $g$ may be not uniformly continuous in $z$. Then, we put forward and prove a comparison theorem and a Levi type theorem on the minimal $L^1$ solutions. A Lebesgue type theorem on $L^1$ solutions is also obtained. Furthermore, we investigate the same problem in the case that $g$ may be discontinuous in $y$. Finally, we prove a general comparison theorem on $L^1$ solutions when $g$ is weakly monotonic in $y$ and uniformly continuous in $z$ as well as it has a stronger sublinear growth in $z$. As a byproduct, we also obtain a general existence and unique theorem on $L^1$ solutions. Our results extend some known works.\vspace{0.1cm}
\end{abstract}

\begin{keyword}
Backward stochastic differential equation \sep Integrable parameters \sep Existence \sep Comparison theorem \sep Levi type theorem \vspace{0.2cm}

\MSC[2010] 60H10
\end{keyword}
\end{frontmatter}\vspace{-0.4cm}


\section{Introduction}

Nonlinear backward stochastic differential equation (BSDE for short) was first introduced in \citep{Par90} by Pardoux and Peng. They established an existence and uniqueness result for solutions to multidimensional BSDEs with square integrable parameters under the Lipschitz assumption of the generator $g$. From then on, BSDEs have been extensively studied, and many applications have been found in mathematical finance, stochastic control, and partial differential equations. Particularly, much effort have been made to relax the Lipschitz hypothesis on $g$, for instance, some results can be found in \cite{Bah02,Bah10,Bri03,Bri06,Bri07,El97,Fan12a,Fan12b,Ham03,Jia08,Kob00,
Lep97,Lep98,Mao95,Par99,Peng91}, most of which dealt with BSDEs with square-integrable parameters.

On the other hand, \citet{Peng97} introduced the notion of $g$-martingales by solutions to BSDEs, which can be viewed, in some sense, as nonlinear martingales. Since the classical theory of martingales is carried in the integrable space, the question of solving a BSDE with only integrable parameters comes up naturally, as has been pointed out in \citet{Bri03}. In recent few years, this question has attracted more and more interests and some important results on it have also been obtained in \citep{Bri02,Bri03,Bri06,Fan07,Fan10,Peng97,Xiao12}.
The objective of this paper is to establish some results in this
direction. We only deal with one-dimensional BSDEs and always assume that both the terminal value $\xi$ and the process $g(t,0,0)$ are only integrable.

In Section 2, we establish the existence for a minimal (maximal) $L^1$ solution of the BSDE when the generator $g$ is stronger continuous in $(y,z)$ and monotonic in $y$ as well as it has a general growth in $y$ and a sublinear growth in $z$ (see Theorem 1 and Remark 3). Particulary, we need neither the Lipschitz continuity assumption nor the H$\ddot{{\rm o}}$lder continuity assumption of $g$ in $z$ required respectively in \citet{Bri03} and \citet{Xiao12}. Hence, Theorem 1 extends the corresponding results (in the one-dimensional case) in two referees quoted before.

In the proof of Theorem 1, we use a localization procedure developed in \citet{Bri06} together with an a prior bound given by the unique $L^1$ solution of a BSDE with a H$\ddot{{\rm o}}$lder continuous generator in $z$. For this purpose, similar to Theorem 4.1 in \citet{Bri07} we construct a sequence of $g_n$ to approach the generator $g$. However, we would like to mention that the sequence $g_n$ is obtained by ``the infinite evolution" made between $g$ and $|z|^{\alpha}$ with $0<\alpha<1$, but not between $g$ and $|z|$ as usual (see Proposition 1 and Remark 2). At the same time, in order to deal with the general growth of $g$ in $y$, we use two stopping time sequences $\{\tau_k\}$ and $\{\sigma_m\}$ different from not only those in Theorem 2.1 of \citet{Fan07} but also those in Theorem 4.1 of \citet{Bri07}. The use of these two stopping time sequences allow us to eliminate the additional continuity assumptions employed in two results quoted before.

Under the same assumptions as in Section 2, we put forward and prove, in Section 3, a comparison theorem and a Levi type theorem on the minimal (maximal) $L^1$ solutions (see Theorems 2-3 and Remark 5). A Lebesgue type theorem on $L^1$ solutions is also obtained in this section (see Theorem 4). We mention that Theorems 3 and 4 improve, in some sense, the main results in \citet{Fan07}.

Section 4 is devoted to the case that the generator $g$ may be discontinuous in $y$. Under the assumptions that $g$ is left-continuous, lower semi-continuous in $y$ and continuous in $z$ as well as it has a linear growth in $y$ and a sublinear growth in $z$, we obtain, as in Sections 2-3, an existence theorem, a comparison theorem and a Levi type theorem on minimal (maximal) $L^1$ solutions (see Theorems 5-7). And we also give a Lebesgue type theorem on $L^1$ solutions (see Theorem 8). Here, we make ``the infinite evolution" between $g$ and $|y|+|z|^{\alpha}$ with $0<\alpha<1$ (see Proposition 2) and use again the localization procedure (see the proof of Theorem 5). We also mention that Theorem 5 extends Theorem 10 in the first version of \citet{Bri06}, and their ideas of the proof are also very different (see Remark 7).

In the last section, by virtue of Theorem 1 in \citet{Fan12a}, we establish a general comparison theorem on $L^1$ solutions when the generator $g$ is weakly monotonic in $y$ and uniformly continuous in $z$ as well as it has a stronger sublinear growth in $z$ (see Theorem 9), which improves the corresponding results in \citet{Fan10} and \citet{Xiao12}. As a byproduct, we also obtain a general existence and unique theorem on $L^1$ solutions when $g$ is stronger continuous in $(y,z)$, monotonic in $y$ and uniformly continuous in $z$ as well as it has a general growth in $y$ and a stronger sublinear growth in $z$ (see Theorem 10), which also extends, in some sense, the corresponding results in \citet{Fan10}, \citet{Xiao12} and \citet{Bri03} (see Remark 11).

Let us close this introduction by giving the notations to be used in all this paper. For the remaining of this paper, let us fix a
nonnegative real number $T>0$ and a positive integer $d$. First of all, $(\Omega,\F,P)$ is a complete probability space carrying a standard $d$-dimensional Brownian motion $(B_t)_{t\geq 0}$. $(\F_t)_{t\geq 0}$ is the natural filtration of the Brownian motion $(B_t)_{t\geq 0}$ augmented by the $P$-null sets of $\F$ and we assume $\F_T=\F$. For every positive integer $n$, we use $| \cdot |$ to denote the norm of Euclidean space $\R^n$. For each real $p>0$, $L^p(\R)$ represents the set of all $\F_T$ -measurable random variable $\xi$ such that $\E[|\xi|^p]<+\infty$, and ${\s}^p$ denotes the set of real-valued, adapted and continuous processes $(Y_t)_{t\in\T}$
such that
$$\|Y\|_{{\s}^p}:=\left( \E[\sup_{t\in\T} |Y_t|^p] \right)
^{1\wedge 1/p}<+\infty. $$
If $p\geq 1$, $\|\cdot\|_{{\s}^p}$ is a norm on ${\s}^p$ and if $p\in (0,1)$,
$(X,X')\longmapsto \|X-X'\|_{{\s}^p}$ defines a distance on
${\s}^p$. Under this metric, ${\s}^p$ is complete. Moreover, let ${\rm M}^p$ denote the set of (equivalent classes of)
$(\F_t)$-progressively measurable, ${\R}^d$-valued processes $\{Z_t,
t\in\T\}$ such that
$$\|Z\|_{{\rm M}^p}:=\left\{ \E\left[\left(\int_0^T |Z_t|^2\
{\rm d}t\right)^{p/2}\right] \right\}^{1\wedge 1/p}<+\infty.
$$
For $p\geq 1$, ${\rm M}^p$ is a Banach space endowed
with this norm and for $p\in (0,1)$, ${\rm M}^p$ is a complete metric space with the resulting distance. We set $\s=\cup_{p>1} {\s}^p$ and let us recall that a continuous process $(Y_t)_{t\in\T}$ belongs to the class (D) if the family $\{ Y_\tau:\tau\in\Sigma_T\}$ is uniformly integrable, where $\Sigma_T$ stands for the set of
all $(\F_t)$-stopping times $\tau$ such that $\tau\leq T$. For a process $Y$ in the class (D), we put
$$\|Y\|_1=\sup\{\E[Y_{\tau}]:\ \tau\in\Sigma_T\}.
$$
The space of $(\F_t)$-progressively measurable continuous processes which belong to the class (D) is complete under this norm.

In this paper, we consider the following one-dimensional BSDE:
\begin{equation}
    y_t=\xi+\int_t^Tg(s,y_s,z_s){\rm d}s-\int_t^Tz_s\cdot {\rm d}B_s,\ \
    t\in\T,
\end{equation}
where $\xi\in \LT$ is called the terminal condition, the random function
$$g(\omega,t,y,z):\Omega\tim \T\tim {\R}\tim {\R}^d \To \R$$
is $(\F_t)$-progressively measurable for each $(y,z)$, called the
generator of BSDE(1). We will sometimes use the notation BSDE$(\xi,g)$ to say that we consider the BSDE whose generator is $g$ and whose terminal condition is $\xi$.

By a solution to BSDE(1) we mean a pair of $(\F_t)$-adapted processes $(y_\cdot,z_\cdot)$ with values in ${\R}\times {\R}^d$ such that $\ps$, $t\mapsto y_t$ is continuous, $t\mapsto z_t$ belongs to ${\rm L}^2(0,T)$, $t\mapsto g(t,y_t,z_t)$ belongs to ${\rm L}^1(0,T)$ and (1) holds true for each $t\in\T$.

If a solution $(y_\cdot,z_\cdot)$ to BSDE(1) satisfies that $y_\cdot$ belongs to the class (D) and $(y_\cdot,z_\cdot)\in \s^\beta\tim\M^\beta$ for any $\beta\in (0,1)$, then it will be called a $L^1$ solution to BSDE(1).

\section{Existence of minimal $L^1$ solutions }

Let us first introduce the following assumptions
on the generator $g$\vspace{0.1cm}:

(H1) $g$ is stronger continuous in $(y,z)$, i.e., $\as$, $\RE\ y,\ z\mapsto g(\omega,t,y,z)$ is
continuous, and $y\mapsto g(\omega,t,y,z)$ is
continuous uniformly with respect to $z$;\vspace{0.1cm}

(H2) $g$ is monotonic in $y$, i.e., there exists
a constant $\mu\geq 0$ such that $\as,\ \RE\ y_1,y_2,z$,
$$(g(\omega,t,y_1,z)-g(\omega,t,y_2,z))(y_1-y_2)\leq\mu
|y_1-y_2|^2;$$

(H3) $g$ has a general growth in $y$, i.e., $\RE\ r\geq 0$, $$\varphi_r(t):=\sup\limits_{|y|\leq r}
|g(\omega,t,y,0)|\in {\rm L}^1(\T\tim
\Omega);$$

(H4) $g$ has a sublinear growth in $z$, i.e., there exist two constants $\lambda>0$ and $\alpha\in (0,1)$ as well as a nonnegative and $(\F_t)$-progressively measurable process $(f_t)_{t\in\T}\in {\rm L}^1(\T\tim\Omega)$ such that $\as$, $\RE\ y, z$,
\begin{equation}
\hspace*{1.5cm}|g(\omega,t,y,z)-g(\omega,t,y,0)|\leq \lambda(f_t(\omega)+|y|+|z|^\alpha);
\end{equation}

(H1') $g$ is continuous in $(y,z)$, i.e., $\as$, $(y,z)\mapsto g(\omega,t,y,z)$ is continuous;\vspace{0.15cm}

(H4') $g$ has a stronger sublinear growth in $z$, i.e., same as (H4) expect that (2) is replaced with
$$|g(\omega,t,y,z)-g(\omega,t,y,0)|\leq \lambda(f_t(\omega)+|y|+|z|)^\alpha;$$

(H4'') $g$ is H$\ddot{{\rm o}}$lder continuous in $z$, uniformly
with respect to $(\omega,t,y)$, i.e., there exist two constants
$\gamma>0$ and $\alpha\in (0,1)$ such that $\as$, $\RE\ y,z_1,z_2$,
$$|g(\omega,t,y,z_1)-g(\omega,t,y,z_2)|\leq
\gamma |z_1-z_2|^{\alpha}.$$

We would like to mention that, to our knowledge, (H2), (H3) together with (H4'), and (H4'') are, respectively, put forward at the first time in \citet{Peng91}, \citet{Bri03} and \citet{Fan10}. But, (H1) and (H4) are new. Note that (H4') will be only used in Section 5.\vspace{0.1cm}

{\bf Remark 1} It is not difficult to see that (H1) is slightly stronger than (H1'). Furthermore, (H2) together with (H4) can imply the following inequality:
$$ g(\omega,t,y,z)\ {\rm sgn}(y)\leq |g(\omega,t,0,0)|+\lambda f_t(\omega)+(\lambda+\mu)|y|+\lambda |z|.$$
Finally, it is easy to verify that (H4'') $\Longrightarrow$ (H4') $\Longrightarrow$ (H4).\vspace{0.2cm}

The main result of this section is as follows.\vspace{0.1cm}

{\bf Theorem 1} (Existence theorem on minimal $L^1$ solutions)\ Let (H1)-(H4) hold true for the generator $g$. Then for each $\xi\in\LT$, BSDE$(\xi,g)$ has a minimal $L^1$ solution $( y_\cdot,z_\cdot)$, i.e, if $(y'_\cdot,z'_\cdot)$ is another $L^1$ solution, then for each $t\in\T$, $$y_t\leq y'_t\ \ \ \ps.$$

{\bf Example 1}\ For each $(\omega,t,y,z)\in \Omega\times\T\times\R\times\R^d$, let
$$g(\omega,t,y,z)=-|B_t(\omega)|\cdot{\text e}^y+(|y|+\sqrt{|z|})\cdot\sin|z|+{1\over \sqrt{t}}1_{t>0}+|B_t(\omega)|^2.\vspace{-0.1cm}$$
It is not hard to check that this $g$ satisfies assumptions (H1)-(H4) with $\mu=1,\lambda=1,\alpha=1/2$ and $f_t(\omega)\equiv 0$. It then follows from Theorem 1 that for each $\xi\in\LT$, BSDE$(\xi,g)$ has a minimal $L^1$ solution.

It should be especially pointed out that this generator $g$ has a general growth in the variable $y$, and it is not uniformly continuous with respect to the variable $z$. So it is of course neither Lipschitz continuous nor H\"{o}lder continuous in $z$. Then, the existence result of $L^1$ solutions to BSDE$(\xi,g)$ with $\xi\in\LT$ can not be obtained by any known results including those in \citep{Bri02,Bri03,Bri06,Fan07,Fan10,Peng97,Xiao12}.\vspace{0.2cm}

Before proving Theorem 1, let us recall the following two lemmas taken from Theorem 1 and Proposition 2 in \citet{Xiao12}.\vspace{0.1cm}

{\bf Lemma 1} (Existence theorem) Let (H1'), (H2)-(H3) and (H4'') hold true for the generator $g$. Then for each $\xi\in\LT$, BSDE$(\xi,g)$ has a unique $L^1$ solution.\vspace{0.2cm}

{\bf Lemma 2} (Comparison theorem) Let $g$ and $g'$ be two generators of BSDEs and one of them satisfies (H2) and (H4''). Let $(y_\cdot,z_\cdot)$ and $(y'_\cdot,z'_\cdot)$ be, respectively, a $L^1$ solution to BSDE$(\xi,g)$ and BSDE$(\xi',g')$. If $\ps,\ \xi\leq \xi'$ and $\as$, for each $(y,z)\in\R\tim\R^d$, $g(t,y,z)\leq g'(t,y,z)$, then for each $t\in \T$,
$$\ y_t\leq y'_t\ \ \ \ps.\vspace{-0.1cm}$$

The following proposition gives a nice approximation of the generator $g$ satisfying (H1)-(H4), which will play an important role in the proof of Theorem 1.\vspace{0.1cm}

{\bf Proposition 1}\ Let (H1)-(H4) hold true for the generator $g$. For each $n\geq 1$ and each $(\omega,t,y,z)\in \Omega\times\T\times\R\times{\R}^d$,
let
\begin{equation}
g_n(\omega,t,y,z)=\inf\limits_{u\in \R^d}
(g(\omega,t,y,u)+(n+\lambda)|u-z|^\alpha),
\end{equation}
where $\lambda$ and $\alpha$ are taken from (H4). Then

(i) For each $n\geq 1$, $g_n(\omega,t,y,z)$ is a mapping from $\Omega\times\T\times\R\times{\R}^d$ into $\R$, and for each $(y,z)$, $g_n(\omega,t,y,z)$ is $(\F_t)$-progressively measurable;

(ii) For each $n\geq 1$ and each $(y,z)$, $\as$, we have
$$g_n(\omega,t,y,z)\leq g_{n+1}(\omega,t,y,z)\leq g(\omega,t,y,z)$$ and
$$|g_n(\omega,t,y,z)-g(\omega,t,y,0)|\leq \lambda(f_t(\omega)+|y|+|z|^\alpha); $$

(iii) For each $n\geq 1$, $g_n$ satisfies (H1'), (H2)-(H3) and (H4'') with $\gamma=n+\lambda$;

(iv) If $(y_n,z_n)\rightarrow(y,z)$ as $n\To \infty$, then when $n\To \infty$, we have
$$ g_n(\omega,t,y_n,z_n)\rightarrow g(\omega,t,y,z)\ \ \ \as.$$

{\bf Proof.} In view of the inequality $|u|^\alpha\leq |u-z|^\alpha+|z|^\alpha$, it follows from (H4) that $\as$, for each $(y,z,u)\in \R^{1+d+d}$,
$$\begin {array}{lll}
g(\omega,t,y,u)+(n+\lambda)|u-z|^\alpha&\geq &
g(\omega,t,y,0)-\lambda (f_t(\omega)+|y|+|u|^{\alpha})+\lambda |u-z|^{\alpha}\\
&\geq & g(\omega,t,y,0)-\lambda (f_t(\omega)+|y|+|z|^\alpha).
\end{array}$$
Thus, $\as$, for each $n\geq 1$ and each $(y,z)\in \R^{1+d}$, $g_n(\omega,t,y,z)$ takes values in $\R$ and
\begin{equation}
g_n(\omega,t,y,z)\geq g(\omega,t,y,0)-\lambda (f_t(\omega)+|y|+|z|^\alpha).
\end{equation}
On the other hand, since the mapping $u\mapsto g(\omega,t,y,u)+(n+\lambda)|u-z|^\alpha$ is continuous and $\Q^d$ is dense in $\R^d$, the infimum in (3) taken over $\R^d$ is equal to the one taken over $\Q^d$. Hence, for each $n\geq 1$ and each $(y,z)\in \R^{1+d}$, $g_n(\omega,t,y,z)$
is $(\F_t)$-progressively measurable for each $(y,z)$. Thus, we get (i).

Furthermore, it follows from (3) and (H4) that for each $n\geq 1$, $\as$,
$$g_n(\omega,t,y,z)\leq g_{n+1}(\omega,t,y,z)\leq g(\omega,t,y,z)\leq g(\omega,t,y,0)+\lambda (f_t(\omega)+|y|+|z|^\alpha).$$
Then (ii) follows from the previous inequality and (4).

In the sequel, we will show (iii). For this, let us recall two basic inequalities:
\begin{equation}
\inf_{x\in D} f(x)-\inf_{x\in D} g(x)\leq \sup_{x\in D}(f(x)-g(x))
\end{equation}
 and
\begin{equation}
|\inf_{x\in D} f(x)-\inf_{x\in D} g(x)|\leq \sup_{x\in
D}|f(x)-g(x)|.\vspace{0.1cm}
\end{equation}
Now, we can prove (iii). First, in view of the inequality $|x|^\alpha-|y|^\alpha\leq |x-y|^\alpha$, it follows from (3) and (6) that
$\as$, for each $(y,z_1,z_2)\in \R\tim \R^d \tim \R^d$,
\begin{equation}
\begin{array}{lll}
\hspace*{-0.5cm}|g_n(\omega,t,y,z_1)-g_n(\omega,t,y,z_2)|&\leq& \Dis\sup_{u\in\R^d}
|(n+\lambda)|u-z_1|^\alpha-(n+\lambda)|u-z_2|^\alpha|\\
&\leq&\Dis (n+\lambda)|z_1-z_2|^\alpha.
\end{array}
\end{equation}
Thus, (H4'') holds true for each $g_n$. Second, by (3) and (6) we can deduce that $\as$, for each $(y,y_0,z)\in \R\times\R\times\R^d$,
$$|g_n(\omega,t,y,z)-g_n(\omega,t,y_0,z)|\leq \sup\limits_{u\in \R^d}|g(\omega,t,y,u)-g(\omega,t,y_0,u)|.$$
Because $\as$, $y\mapsto g(\omega,t,y,z)$ is continuous uniformly with respect to $z$ by (H1), from the previous inequality we know that $\as$, for each $z\in \R^d$, $y\mapsto g_n(\omega,t,y,z)$ is continuous. On the other hand, (7) means that $\as$, $z\mapsto g_n(\omega,t,y,z)$ is uniformly continuous uniformly with respect to $y$. Hence, we can conclude that $\as$, $(y,z)\mapsto g_n(\omega,t,y,z)$ is continuous, that is, (H1') holds true for each $g_n$. Third, in view of (H2), it follows from (3) and (5) that, $\as$, for each $(y_1,y_2)\in\R\tim\R$ with $y_1\geq y_2$,
$$\begin{array}{lll}
\Dis &&(y_1-y_2)(g_n(\omega,t,y_1,z)-g_n(\omega,t,y_2,z))\\
 &\leq&
\Dis\sup_{u\in \R^d}\{(y_1-y_2)(g_n(\omega,t,y_1,u)-
g_n(\omega,t,y_2,u))\}\leq \mu |y_1-y_2|^2.
\end{array}
$$
Furthermore, if $y_1<y_2$, then by exchanging the position of $y_1$ and $y_2$ we know that the above inequality holds also true. Therefore, (H2) is also true for $g_n$. At last, it follows from (ii) that, $\as$, for each $n\geq 1$ and each $y\in \R$,
$$|g_n(\omega,t,y,0)|\leq |g(\omega,t,y,0)|+ \lambda(f_t(\omega)+|y|).$$
Thus, $g_n$ also satisfies (H3) since $g$ satisfies it. (iii) is then proved.

Finally, we prove (iv). Assume that
$(y_n,z_n)\longrightarrow (y,z)$ as $n\To\infty$. By (3) and (H4) we can take a sequence $(v_n)_{n\geq 1}$ such that $\as$,
\begin{equation}
\begin{array}{lll}
&& g_n(\omega,t,y_n,z_n)\\
&\geq& g(\omega,t,y_n,v_n)+(n+\lambda)|z_n-v_n|^\alpha-{1\over n}\\ &\geq& g(\omega,t,y_n,0)-\lambda(f_t+|y_n|+|v_n|^\alpha)+
(n+\lambda)|z_n-v_n|^\alpha-{1\over n}\\ &\geq& g(\omega,t,y_n,0)-\lambda(f_t+|y_n|+|z_n|^\alpha)
+(n+\lambda-1)|z_n-v_n|^\alpha-{1\over n}.  \end{array}
\end{equation}
Furthermore, it follows from (ii) that
$$g_n(\omega,t,y_n,z_n)\leq g(\omega,t,y_n,0)+\lambda(f_t(\omega)+|y_n|+|z_n|^\alpha).$$
Thus, we have
$$(n+\lambda-1)|z_n-v_n|^\alpha\leq 2\lambda(f_t(\omega)+|y_n|+|z_n|^\alpha)+\frac{1}{n},$$
and then
$$\limsup_{n\to\infty} n|z_n-v_n|^\alpha<+\infty.$$
Therefore
$$\lim\limits_{n\to\infty}v_n=z.\vspace{0.2cm}$$
Then, in view of Remark 1, it follows from (8) and (H1) that $\as$,
$$\liminf_{n\to\infty}g_n(\omega,t,y_n,z_n)\geq\liminf_{n\to\infty} g(\omega,t,y_n,v_n)=g(\omega,t,y,z).$$
On the other hand, it follows from (ii) and (H1), in view of Remark 1, that $\as$,
$$\limsup_{n\to\infty} g_n(\omega,t,y_n,z_n)\leq\limsup_{n\to\infty} g(\omega,t,y_n,z_n)=g(\omega,t,y,z).$$
Hence, we have (iv), and Proposition 1 is proved.\vspace{0.2cm}\hfill $\Box$

{\bf Remark 2} Similar argument to Proposition 1 yields that if we replace (3) with
$$
g^n(\omega,t,y,z)=\sup\limits_{u\in \R^d}
(g(\omega,t,y,u)-(n+\lambda)|u-z|^\alpha),
$$
then the conclusions of Proposition 1 hold also true for $g^n$, except that $g^n$ is non-increasing with respect to $n$ and bigger than $g$.\vspace{0.2cm}

Now we can turn to the proof of Theorem 1.\vspace{0.1cm}

{\bf Proof of Theorem 1.} Suppose that $\xi\in\LT$ and that (H1)-(H4) hold for the generator $g$. For each $n\in\N$ and each $(\omega,t,y,z)\in \Omega\tim \T\tim \R\tim \R^d$, let $g_n(\omega,t,y,z)$ be defined in (3) and $h(\omega,t,y,z)$ be defined as follows
$$h(\omega,t,y,z):=g(\omega,t,y,0)+\lambda(f_t(\omega)+|y|+
|z|^{\alpha}),\ \ \RE\ \omega,t,y,z.$$
Since $g$ satisfies (H1)-(H4), in view of (ii) and (iii) in Proposition 1 and Remark 1, we know that both $h$ and $g_n$ satisfy (H1'), (H2)-(H3) and (H4''), and $\as$, for each $n\geq 1$ and each $(y,z)\in \R\tim\R^d$, we have
$$g_n(\omega,t,y,z)\leq g_{n+1}(\omega,t,y,z)\leq g(\omega,t,y,z)\leq h(\omega,t,y,z).$$
Then, it follows from Lemma 1 that for each $n\geq 1$, both BSDE$(\xi,g_n)$ and BSDE$(\xi,h)$ have unique $L^1$ solutions, denoted, respectively,  by $(y^n_t,z^n_t)_{t\in \T}$ and $(\tilde{y}_t,\tilde{z}_t)_{t\in \T}$ for notational convenience. Furthermore, by Lemma 2 we also know that for each $t\in \T$,
$$
  y^1_t\leq y^n_t\leq y^{n+1}_t\leq \tilde{y}_t\ \ \ \ps.
$$
We define $y_\cdot=\Lim y^n_\cdot$, then
\begin{equation}
\RE\ t\in[0,T],\ \ |y_t|\leq \sup\limits_{n\geq 1}|y^n_t|\leq |y^1_t|+|\tilde{y}_t|\ \ \ \ps.\vspace{-0.1cm}
\end{equation}

In the sequel, we will use a similar localization procedure as in \citet{Bri06}. For each $k\in\N$, let us introduce the following stopping time:
$$\tau_k=\inf\{t\in\T:|y^1_t|+|\tilde{y}_t|\geq k\}\wedge T.$$
For fixed $k$ and $m\in\N$, define also the
following stopping time:
$$\sigma_m=\inf\{t\in\T:\int_0^t (\varphi_k(u)+f_u)\ {\rm d}u\geq m\}\wedge T,$$
where $\varphi_k(u)$ and $f_u$ are defined in assumptions (H3) and (H4) respectively. Then $(y^n_{k,m}(t),z^n_{k,m}(t)):=(y^n_{t\wedge (\tau_k \wedge\sigma_m)},z^n_t 1_{t\leq (\tau_k \wedge \sigma_m)})$
solves the following BSDE:
$$
  y_{k,m}^n(t)=\xi^n_{k,m}+\int_t^T 1_{u\leq (\tau_k \wedge \sigma_m)} g_n(u,y_{k,m}^n(u),z_{k,m}^n(u))
  {\rm d}u-\int_t^Tz_{k,m}^n(u)\cdot {\rm d}B_u,
$$
where $\xi^n_{k,m}:=y^n_{k,m}(T)=y^n_{\tau_k \wedge
\sigma_m}$.\vspace{0.1cm}

For each pair of $k$ and $m$, it is very important to observe that, $y_{k,m}^n$ is nondecreasing in $n$ by construction. Further, it follows from the definition of $\tau_k$ and the inequality (9) that
\begin{equation}
\sup\limits_{n\geq 1}\sup\limits_{t\in[0,T]}\|y_{k,m}^n(t)
\|_\infty\leq k.
\end{equation}
Then, in view of this inequality, Remark 1 and the definitions of $\tau_k$  and $\sigma_m$ again, it follows from Lemma 3.1 in \citet{Bri03} that
$$z_{k,m}^n\in \M^2,\ \ \RE\ n\in\N.$$
Thus, if $\rho_k(y)=yk/\max(|y|,k)$, we have
\begin{equation}
y_{k,m}^n(t)=\xi^n_{k,m}+\int_t^T 1_{u\leq (\tau_k \wedge \sigma_m)} g_n(u,\rho_k(y_{k,m}^n(u)),z_{k,m}^n(u)){\rm d}u-\int_t^Tz_{k,m}^n(u)\cdot {\rm d}B_u.
\end{equation}
Moreover, by (ii) of Proposition 1 and (H3) we have
$$
\begin{array}{lll}
|1_{u\leq (\tau_k \wedge \sigma_m)}g_n(u,\rho_k(y),z)|
&\leq & 1_{u\leq (\tau_k \wedge \sigma_m)}\left[|g(u,\rho_k(y),0)|
+\lambda(f_u+|\rho_k(y)|+|z|^\alpha)\right]\\
&\leq & 1_{u\leq (\tau_k \wedge \sigma_m)}\left[\varphi_k(u)+\lambda(f_u+k+1+|z|)
\right].\\
&\leq & (1+\lambda)1_{u\leq (\tau_k \wedge \sigma_m)}(\varphi_k(u)+f_u)+\lambda(k+1)+\lambda |z|.
\end{array}$$
It then follows from the definition of $\sigma_m$ that
$$\int_0^T 1_{u\leq (\tau_k \wedge \sigma_m)}(\varphi_k(u)+f_u)\ {\rm d}u \leq m. $$
In view of the previous two inequalities, (10), (iv) of Proposition 1 and the fact that $y_{k,m}^n$ is nondecreasing in $n$, arguing as in the proof of Theorem 1 in \citet{Lep97} (see pages 427-429), we can take the limit with respect to $n$ ($k$ and $m$ being fixed) in (11) in the space $\s^2\tim\M^2$, where the only change need to be made is that we have to use Lebesgue's dominated convergence theorem in stead of H\"{o}lder inequality in order to show the convergence of $z^n_{k,m}$ in ${\rm M}^2$. In particular, setting $y_{k,m}=\Lim y_{k,m}^n$, we know that $y_{k,m}$ is continuous and that there exists a process $z_{k,m}\in \M^2$ such that $\Lim
z^n_{k,m}=z_{k,m}$ in $\M^2$ and $(y_{k,m},z_{k,m})$ solves
$$
y_{k,m}(t)=\xi_{k,m}+\int_t^T 1_{u\leq (\tau_k \wedge \sigma_m)} g(u,\rho_k(y_{k,m}(u)),z_{k,m}(u)){\rm d}u-\int_t^Tz_{k,m}(u)\cdot {\rm d}B_u,\ t\in\T,
$$
where $\xi_{k,m}:=\Lim y^n_{\tau_k \wedge \sigma_m}$.

Since
$|y_{k,m}(t)|\leq k$, the above equation can be rewritten as
\begin{equation}
y_{k,m}(t)=\xi_{k,m}+\int_t^T 1_{u\leq (\tau_k \wedge \sigma_m)} g(u,y_{k,m}(u),z_{k,m}(u)){\rm d}u-\int_t^Tz_{k,m}(u)\cdot {\rm d}B_u.
\end{equation}
But $\sigma_m\leq \sigma_{m+1}$, $\tau_k\leq \tau_{k+1}$, then we get, using the definition of $y_{k,m},z_{k,m}$ and $y_\cdot$,
$$
\begin{array}{l}
y_{t\wedge\tau_k \wedge \sigma_m}=y_{k+1,m+1}(t\wedge
\tau_k\wedge\sigma_m)=y_{k,m}(t)=\Lim y_{t\wedge\tau_k \wedge \sigma_m}^n,\\
z_{k+1,m+1}(t)1_{t\leq (\tau_k \wedge \sigma_m)}=z_{k,m}(t)=\Lim
z^n_t1_{t\leq (\tau_k \wedge \sigma_m)}.
\end{array}$$
It follows from (H3), (H4) and the definitions of $\tau_k$ and $\sigma_m$ that $\sigma_m\To T$ as $m\To \infty$ for each fixed $k$ and $\tau_k\To T$ as $k\To \infty$, and
thus since all of $y_{k,m}$ are continuous processes we deduce that $y_\cdot$ is continuous on $[0,T]$. Then we define $z_\cdot$ on $(0,T)$ by setting
$$z_t=z_{k,m}(t),\ \ {\rm if}\ t\in(0,\tau_k \wedge \sigma_m),\vspace{0.2cm}$$
so that $z_t1_{t\leq (\tau_k \wedge \sigma_m)}=z_{k,m}(t)1_{t\leq (\tau_k \wedge \sigma_m)}=z_{k,m}(t)$ and (12) can be rewritten as
\begin{equation}
y_{t\wedge\tau_k \wedge \sigma_m}=y_{\tau_k \wedge \sigma_m}  +\int_{t\wedge\tau_k \wedge \sigma_m}^{\tau_k \wedge \sigma_m} g(u,y_u,z_u) {\rm d}u-\int_{t\wedge\tau_k \wedge \sigma_m}^{\tau_k \wedge  \sigma_m}z_u\cdot {\rm d}B_u.
\end{equation}

Furthermore, we have
$$
\begin{array}{lll}
&&\Dis P\left(\int_0^T|z_u|^2{\rm d}u=\infty\right)\\
&=&\Dis P \left(\int_0^T|z_u|^2{\rm d}u=\infty,\tau_k \wedge
\sigma_m=T\right)+ P \left(\int_0^T|z_u|^2{\rm
d}u=\infty,\tau_k
\wedge \sigma_m<T\right)\\
&\leq & \Dis P \left(\int_0^{\tau_k \wedge
\sigma_m}|z_{k,m}(u)|^2{\rm d}u=\infty\right)+P\left(\tau_k \wedge \sigma_m<T\right),
\end{array}$$
and we deduce, in view of the fact that $\sigma_m\To T$ as $m\To \infty$ for each fixed $k$, and $\tau_k\To T$ as $k\To \infty$, that
$$\int_0^T|z_u|^2{\rm d}u<\infty\ \ \ \ps.$$
Let $m\To \infty$ for fixed $k$ in (13), and then let $k\To\infty$, we deduce that $(y_\cdot,z_\cdot)$ is a solution of BSDE$(\xi,g)$. By (9) we know that $y_\cdot$ belongs to the class (D) and the space $\s^\beta$ for each $\beta\in (0,1)$. Furthermore, in view of Remark 1, by Lemma 3.1 in \citet{Bri03} we also know that $z_\cdot$ belongs to the space $\M^\beta$ for each $\beta\in (0,1)$. Consequently, $(y_\cdot,z_\cdot)$ is a $L^1$ solution of BSDE$(\xi,g)$.

Finally, we prove that $(y_\cdot,z_\cdot)$ is also a minimal $L^1$ solution of BSDE$(\xi,g)$. Let $(y'_\cdot,z'_\cdot)$ be another $L^1$ solution. Note that for each $n\geq 1$, $g_n$ satisfies (H2) and (H4''), and it is smaller than $g$. It follows from Lemma 2 that for each $t\in\T$ and each $n\geq 1$, $y^n_t\leq y'_t\ \ps$. Since $y_\cdot=\Lim y^n_\cdot$, we can obtain that for each $t\in\T$, $y_t\leq y'_t\ \ps$. The proof of Theorem 1 is then completed.\hfill $\Box$\vspace{0.2cm}

{\bf Remark 3}\ In the proof of Theorem 1, replace $g_n$ with $g^n$ defined in Remark 2, and $h$ with the following function:
$$\tilde{h}(\omega,t,y,z):=g(\omega,t,y,0)-\lambda(f_t(\omega)
+|y|+|z|^{\alpha}),$$
and let $(Y^n_t,Z^n_t)_{t\in\T}$ be the unique $L^1$ solution of BSDE$(\xi,g^n)$ by virtue of Remark 2 and Lemma 1. Then, using the similar procedure as that in the proof of Theorem 1, we can deduce that the limit process $(Y_\cdot,Z_\cdot)$ of the sequence $(Y^n_\cdot,Z^n_\cdot)$ is a maximal $L^1$ solution of BSDE$(\xi,g)$ when $g$ satisfies (H1)-(H4), i.e, if $(Y'_\cdot,Z'_\cdot)$ is another $L^1$ solution, then for each $t\in\T$, $Y'_t\leq Y_t\ \ps$.

\section{Comparison theorem and Levi type theorem on minimal $L^1$ solutions}

In this section, we will put forward and prove a comparison theorem (Theorem 2) and a Levi type theorem (Theorem 3) on the minimal $L^1$ solution of BSDE(1) under (H1)-(H4). A Lebesgue type theorem (Theorem 4) on $L^1$ solutions is also obtained in this section.

We mention that Theorems 3 and 4 improve, in some sense, the corresponding results of \citet{Fan07} since the additional continuity assumption (H) and the Lipschitz continuity assumption of $g$ in $z$ employed in \cite{Fan07} are moved away in Theorems 3 and 4, and the assumption (H4') used in \cite{Fan07} is also weakened to (H4) here.\vspace{0.1cm}

{\bf Theorem 2} (Comparison theorem on the minimal $L^1$ solution)\ Assume that $\xi,\xi'\in \LT$ and that both $g$ and $g'$ satisfy (H1)-(H4). Let $(y^{g}_t(\xi),z^{g}_t(\xi))_{t\in \T}$ and $(y^{g'}_t(\xi'),z^{g'}_t(\xi'))_{t\in \T}$ be, respectively, the minimal $L^1$ solution to BSDE$(\xi,g)$ and BSDE$(\xi',g')$ by Theorem 1. If $\ps,\ \xi\leq \xi'$ and $\as$, for each $(y,z)\in \R\tim \R^d$, $g(t,y,z)\leq g'(t,y,z)$, then for each $t\in\T$, we have
$$y^{g}_t(\xi)\leq y^{g'}_t(\xi')\ \ \ \ps.$$

{\bf Proof.}\ Let $g_n$ be defined in (3). By (iii) of Proposition 1 and the proof procedure of Theorem 1 we know that for each $n\geq 1$, $g_n$ satisfies (H1'), (H2)-(H3) and (H4''), and for each $t\in\T$, \begin{equation}
\Lim y^{g_n}_t(\xi)=y^{g}_t(\xi)\ \ \ \ps,
\end{equation}
where $(y^{g_n}_t(\xi),z^{g_n}_t(\xi))_{t\in\T}$ is the unique $L^1$ solution of BSDE$(\xi,g_n)$.

On the other hand, in view of the assumptions of Theorem 2, by (ii) of Proposition 1 we also know that $\ps, \xi\leq \xi'$ and $\as$, for each $n\geq 1$ and $(y,z)\in \R\tim\R^d$,
$$g_n(t,y,z)\leq g(t,y,z)\leq g'(t,y,z).\vspace{0.1cm}$$
Then, noticing that $g_n$ satisfies (H2) and (H4''), by Lemma 2 we get that for each $n\geq 1$ and $t\in \T$,
\begin{equation}
y^{g_n}_t(\xi)\leq y^{g'}_t(\xi')\ \ \ \ps.\vspace{0.1cm}
\end{equation}
Thus, the conclusion of Theorem 2 follows from (14) and (15).\hfill $\Box$\vspace{0.2cm}

{\bf Theorem 3} \ (Levi type theorem on the minimal $L^1$ solution)\ Assume that $\xi_n,\xi\in \LT$ for each $n\geq 1$ and that $g$ satisfies (H1)-(H4). Let $(y^{g}_t(\xi_n),z^{g}_t(\xi_n))_{t\in \T}$ and $(y^{g}_t(\xi),z^{g}_t(\xi))_{t\in \T}$ be, respectively, the minimal $L^1$ solution of BSDE$(\xi_n,g)$ and BSDE$(\xi,g)$ by Theorem 1. If $\ps,\ \xi_n\uparrow\xi$, then for each $t\in\T$,
$$\Lim\uparrow y^{g}_t(\xi_n)=y^{g}_t(\xi)\ \ \ \ps.$$

{\bf Proof.}\ In view of $\xi_n\uparrow\xi\ \ps$, it follows from Theorem 2 that for each $t\in \T$ and $n\geq 1$,
\begin{equation}
y^{g}_t(\xi_1)\leq y^{g}_t(\xi_n)\leq y^{g}_t(\xi_{n+1})\leq y^{g}_t(\xi)\ \ \ \ps.
\end{equation}
We define $y_\cdot=\Lim y^{g}_\cdot(\xi_n)$, then
$$
\RE\ t\in[0,T],\ \ |y_t|\leq \sup\limits_{n\geq 1}|y^{g}_t(\xi_n)|\leq |y^{g}_t(\xi_1)|+|y^{g}_t(\xi)|\ \ \ \ps.
$$
Thus, for each $k\in\N$, we introduce the following stopping time:
$$\tau_k=\inf\{t\in\T:|y^{g}_t(\xi_1)|+|y^{g}_t(\xi)|+\int_0^t |g(u,0,0)|\ {\rm d}u\geq k\}\wedge T,$$
and for fixed $k$ and $m\in\N$, let the
stopping time $\sigma_m$ be defined in the proof of Theorem 1. Then
$$(y^n_{k,m}(t),z^n_{k,m}(t)):=(y^g_{t\wedge (\tau_k \wedge\sigma_m)}(\xi_n),z^g_t(\xi_n) 1_{t\leq (\tau_k \wedge \sigma_m)})\vspace{0.2cm}$$
solves the following BSDE
$$
  y_{k,m}^n(t)=\xi^n_{k,m}+\int_t^T 1_{u\leq (\tau_k \wedge \sigma_m)} g(u,y_{k,m}^n(u),z_{k,m}^n(u))
  {\rm d}u-\int_t^Tz_{k,m}^n(u)\cdot {\rm d}B_u,
$$
where $\xi^n_{k,m}:=y^n_{k,m}(T)=y^g_{\tau_k \wedge
\sigma_m}(\xi_n)$.

In the sequel, arguing as in the proof of Theorem 1, we can deduce that there exists a process $z_\cdot$ such that $(y_\cdot,z_\cdot)$ is a $L^1$ solution of BSDE$(\xi,g)$. Furthermore, in view of (16), the definition of $y_\cdot$ and the fact that $(y^{g}_\cdot(\xi),z^{g}_\cdot(\xi))$ be the minimal $L^1$ solution of BSDE$(\xi,g)$, we know that $$(y_\cdot,z_\cdot)=(y^{g}_\cdot(\xi),z^{g}_\cdot(\xi)), $$ from which the conclusion of Theorem 3 follows immediately.\hfill $\Box$\vspace{0.2cm}

{\bf Remark 4}\  If the condition ``$\xi_n\uparrow \xi$" in Theorem 3 is replaced with ``$\xi_n\downarrow \xi$", then the sign ``$\leq$" in (16) will change to ``$\geq$", and the $(y_\cdot,z_\cdot)$ in the proof of Theorem 3 is still a $L^1$ solution of BSDE$(\xi,g)$, but
it is uncertain whether it is the minimal one or not, so the conclusion of Theorem 3 does not hold in general. However, if we further assume that the $L^1$ solution of BSDE$(\xi,g)$ is unique, then the conclusion will hold.\vspace{0.2cm}

{\bf Remark 5}\ Using the similar arguments as in Theorems 2-3, in view of Remark 3, we can prove that, in Theorems 2-3, if we replace the minimal $L^1$ solution with the maximal $L^1$ solution, and ``$\xi_n\uparrow \xi$" with ``$\xi_n\downarrow \xi$", then the conclusions hold also true.\vspace{0.2cm}

If the $L^1$ solution of BSDE$(\xi,g)$ is unique, we have the following Lebesgue type theorem on the $L^1$ solution.\vspace{0.1cm}

{\bf Theorem 4}\ (Lebesgue type theorem on the $L^1$ solution)\ Assume that $\xi_n,\xi\in \LT$ for each $n\geq 1$ and that $g$ satisfies (H1)-(H4). Assume further that BSDE$(\xi,g)$ has a unique $L^1$ solution $(y^{g}_t(\xi),z^{g}_t(\xi))_{t\in \T}$. Let $(y^{g}_t(\xi_n),z^{g}_t(\xi_n))_{t\in \T}$ be any of $L^1$ solutions of BSDE$(\xi_n,g)$ by Theorem 1 and Remark 3. If $\ps,\ \xi_n\To\xi$ as $n\To\infty$ and $\ps,\ |\xi_n|\leq \eta$ with $\E[|\eta|]<+\infty$, then for each $t\in\T$,
$$\Lim y^{g}_t(\xi_n)=y^{g}_t(\xi)\ \ \ \ps.$$

{\bf Proof.}\ Let
$$\bar{\xi}_n:=\sup\limits_{k\geq n}\xi_k\ \ {\rm and}\ \
\underline{\xi}_n:=\inf\limits_{k\geq n}\xi_k.\vspace{0.2cm}$$
Then, both $\bar{\xi}_n$ and $\underline{\xi}_n$ belongs to $\LT$ since $\ps,\ |\xi_n|\leq \eta$ with $\E[|\eta|]<+\infty$. And, since $\xi_n\To\xi\ \ps$ as $n\To\infty$, we have, $\ps$,
\begin{equation}
\underline{\xi}_n\leq \xi_n\leq\bar{\xi}_n,\ \ \bar{\xi}_n\downarrow \xi\ \ {\rm and}\ \ \underline{\xi}_n\uparrow \xi.
\end{equation}
In view of Theorem 1, we can let  $$(\underline{y}^{g}_t(\xi_n),\underline{z}^{g}_t(\xi_n))_{t\in \T},\ \ (\underline{y}^{g}_t(\underline{\xi}_n),
\underline{z}^{g}_t(\underline{\xi}_n))_{t\in \T}\ \ {\rm and}\ \ (\underline{y}^{g}_t(\bar{\xi}_n),\underline{z}^{g}_t
(\bar{\xi}_n))_{t\in \T},$$
respectively, be the minimal $L^1$ solution of BSDE$(\xi_n,g)$, BSDE$(\underline{\xi}_n,g)$ and BSDE$(\bar{\xi}_n,g)$. Then, in view of (17) and the fact that $(y^{g}_t(\xi),z^{g}_t(\xi))_{t\in \T}$ is the unique $L^1$ solution of BSDE$(\xi,g)$, by Theorem 3, Theorem 2 and Remark 4 we can deduce that
$$y^{g}_t(\xi)=\Lim\uparrow\underline{y}^{g}_t(\underline{\xi}_n)
\leq \llim \underline{y}^{g}_t(\xi_n)\leq \slim \underline{y}^{g}_t(\xi_n)\leq \Lim\downarrow
\underline{y}^{g}_t(\bar{\xi}_n)=y^{g}_t(\xi),$$
which means that for each $t\in\T$,
$$\Lim \underline{y}^{g}_t(\xi_n)=y^{g}_t(\xi)\ \ \ \ps.$$
In the same way, in view of Remark 5, we can also prove that for each $t\in\T$,
$$\Lim \bar{y}^{g}_t(\xi_n)=y^{g}_t(\xi)\ \ \ \ps,$$
where $(\bar{y}^{g}_t(\xi_n),\bar{z}^{g}_t(\xi_n))$ represents the maximal $L^1$ solution of BSDE$(\xi_n,g)$. Thus, the conclusion of Theorem 4 follows from the above last two identities.\hfill $\Box$

\section{The case that $g$ may be discontinuous in $y$ }

Let us further introduce the following assumptions:\vspace{0.1cm}

(H5) $g$ has a linear growth in $y$ and a sublinear growth in $z$, i.e., there exists two constants $C> 0$, $\alpha\in (0,1)$ and a non-negative $(\F_t)$-progressively measurable stochastic process $(f_t)_{t\in \T}\in {\rm L}^1(\T\tim\Omega)$ such that $\as$,
$$|g(\omega,t,y,z)|\leq f_t(\omega)+
C(|y|+|z|^\alpha),\ \ \RE\ (y,z)\in \R^{1+d}.$$

(H1a) $g$ is left-continuous and lower semi-continuous in $y$, and continuous in $z$, i.e., $\as$, for each $(y_0,z_0)\in \R^{1+d}$, we have
\begin{equation}
\lim\limits_{(y,z)\To (y_0^-, z_0)}g(\omega,t,y,z)=g(\omega,t,y_0,z_0)
\end{equation}
and
\begin{equation}
\liminf\limits_{(y,z)\To (y_0^+,z_0)}g(\omega,t,y,z)\geq g(\omega,t,y_0,z_0).\vspace{0.1cm}
\end{equation}

(H1b) $g$ is right-continuous and upper semi-continuous in $y$, and continuous in $z$, i.e., $\as$, for each $(y_0,z_0)\in \R^{1+d}$, we have
$$
\lim\limits_{(y,z)\To (y_0^+, z_0)}g(\omega,t,y,z)=g(\omega,t,y_0,z_0)$$
and
$$\limsup\limits_{(y,z)\To (y_0^-,z_0)}g(\omega,t,y,z)\leq g(\omega,t,y_0,z_0).\vspace{0.1cm}
$$

{\bf Remark 6}\ Note that (H1a) and (H1b) are taken from \citet{Fan12b}, where the $L^2$ solutions to BSDEs are investigated when $g$ satisfies (H1a) (or (H1b)) and (H5) with $\alpha=1$. It is clear that (H1a)$+$(H1b) $\Leftrightarrow$ (H1'). If (H1a) (resp. (H1b)) holds for $g$, then, $\as$, for each $(y_0,z_0)\in \R\times \R^d$,
$$\liminf_{(y,z)\To (y_0,z_0)}g(\omega,t,y,z)\geq g(\omega,t,y_0,z_0)$$
$$({\rm resp.}\limsup_{(y,z)\To (y_0,z_0)}g(\omega,t,y,z)\leq g(\omega,t,y_0,z_0)).
$$
But $g$ may be discontinuous in $y$ when (H1a) or (H1b) holds true for it. In addition, by virtue of the knowledge of mathematical analysis it is not hard to conclude that if $\as$, for each $z\in \R^d$, $y\mapsto g(\omega,t,y,z)$ is left-continuous (resp. right-continuous) and nondecreasing, and $\as$, for each $y\in \R$, $z\mapsto g(\omega,t,y,z)$ is also continuous, then $g$ must satisfy (H1a) (resp. (H1b)) (see Section 3 in \citet{Fan12b} for more details).\vspace{0.2cm}

The following Theorem 5 establishes an existence result on minimal $L^1$ solutions of BSDEs with discontinuous generators in $y$, which is one of the main results of this section.\vspace{0.1cm}

{\bf Theorem 5} (Existence theorem on the minimal (resp. maximal) $L^1$ solution)\ Assume that the generator $g$ satisfies (H1a) (resp. (H1b)) and (H5). Then for each $\xi\in \LT$, BSDE$(\xi,g)$ has a minimal (resp. maximal) $L^1$ solution $(y_\cdot,z_\cdot)$.\vspace{0.2cm}

{\bf Remark 7}\ A similar result to Theorem 5 was obtained in Theorem 10 of the first version of \citet{Bri06}, where the generator $g$ is continuous in $(y,z)$ and the $f_t(\omega)$ in (H5) is a constant. In addition, it should be mentioned that the $L^1$ solution of BSDE$(\xi,g)$ constructed by them is not necessarily the minimal or maximal one. Hence, Theorem 5 extends this known result.

At the same time, the basic idea developed in Theorem 10 of the first version of \citet{Bri06} is to approach the $L^1$ solution of BSDE$(\xi,g)$ by virtue of a $L^2$ solution sequence of BSDE$(\xi_{n,p},g)$, where $\xi_{n,p}:=\xi^+\wedge n-\xi^{-}\wedge p$. Compared with it, a very different idea will be employed to prove our Theorem 5. More specifically, we will approach the $L^1$ solution of BSDE$(\xi,g)$ by virtue of a $L^1$ solution sequence of BSDE$(\xi,g_n)$, where the sequence $g_n$ is obtained by ``the infinite evolution" made between $g$ and $|y|+|z|^{\alpha}$. \vspace{0.2cm}

{\bf Example 2}\ For each $(\omega,t,y,z)\in \Omega\times\T\times\R\times\R^d$, let
$$g(\omega,t,y,z)=1_{y\leq 0}\sin y+1_{y>0}\cos y+[|y|+\ln(1+|z|)]\cdot\sin(y^2|z|^3)+B_t(\omega).\vspace{-0.1cm}$$
It is clear that $g$ is discontinuous in $y$ and not uniformly continuous in $z$. It is also easy to verify that $g$ satisfies (H1a) and (H5) with $C=1$ and any $\alpha\in (0,1)$. It then follows from Theorem 5 that for each $\xi\in\LT$, BSDE$(\xi,g)$ has a minimal $L^1$ solution. Note that this conclusion can not be obtained by any existing result.\vspace{0.2cm}

In the proof of Theorem 5, the following Proposition 2 will play an important role, which gives a nice approximation of $g$ satisfying (H1a) and (H5).\vspace{0.1cm}

{\bf Proposition 2}\ Let (H1a) and (H5) hold true for the generator $g$. For each $n\geq 1$ and each $(\omega,t,y,z)\in \Omega\times\T\times\R\times{\R}^d$,
let
\begin{equation}
g_n(\omega,t,y,z):=\inf\limits_{(u,v)\in {\Q}^{1+d}}
\{g(\omega,t,u,v)+nC(|y-u|+|z-v|^\alpha)\}.
\end{equation}
where $C$ and $\alpha$ are taken from (H5). Then

(i) For each $n\geq 1$, $g_n(\omega,t,y,z)$ is a mapping from $\Omega\times\T\times\R\times{\R}^d$ into $\R$, and for each $(y,z)$, $g_n(\omega,t,y,z)$ is $(\F_t)$-progressively measurable;

(ii) For each $n\geq 1$ and each $(y,z)$, $\as$, we have
$$g_n(\omega,t,y,z)\leq g_{n+1}(\omega,t,y,z)\leq g(\omega,t,y,z)$$ and
$$|g_n(\omega,t,y,z)|\leq f_t(\omega)+C(|y|+|z|^\alpha); $$

(iii) For each $ y_1,y_2,z_1,z_2$, $\as$, we
have
$$|g_n(\omega,t,y_1,z_1)-g_n(\omega,t,y_2,z_2)|\leq
nC(|y_1-y_2|+|z_1-z_2|^\alpha);$$

(iv) If $(y_n,z_n)\rightarrow (y_0^-,z_0)$ as $n\To \infty$, then
$$\Lim g_n(\omega,t,y_n,z_n)=g(\omega,t,y_0,z_0)\ \ \ \as.
$$

{\bf Proof.} \ In view of the inequalities $|v|^\alpha\leq |v-z|^\alpha+|z|^\alpha$ and $|u|\leq |y-u|+|y|$, it follows from (20) and (H5) that for each $n\geq 1$, $\as$, for each $(y,z)\in \R^{1+d}$,
$$\begin{array}{lll}
g_n(\omega,t,y,z) &\geq & \Dis \inf_{(u,v)\in \R^{1+d}}\{-f_t(\omega)-C|u|-C|v|^\alpha+C(|y-u|+
|z-v|^{\alpha})\}\\
&\geq &\Dis -f_t(\omega)-C(|y|+|z|^\alpha)
\end{array}$$
and
$$g_n(\omega,t,y,z)\leq g(\omega,t,y,z)\leq f_t(\omega)+C(|y|+|z|^\alpha).\vspace{0.2cm}$$
Thus, (i) and (ii) follows immediately by (20). Furthermore, (iii) follows from (20), (6) and the basic inequality $|x_1|^\alpha-|x_2|^\alpha\leq |x_1-x_2|^\alpha$.

Hence, it suffices to show (iv). Indeed, assume that
$(y_n,z_n)\rightarrow (y_0^-,z_0)$ as $n\To \infty$.
In view of the inequalities $|v_n|^\alpha\leq |z_n-v_n|^\alpha+|z_n|^\alpha$ and $|u_n|\leq |y_n-u_n|+|y_n|$, from (20) and (H5) we can take a sequence $(u_n,v_n)$ such that $\as$,
\begin{equation}
\begin{array}{lll}
&& g_n(\omega,t,y_n,z_n)\\
&\geq & \Dis
g(\omega,t,u_n,v_n)+nC(|y_n-u_n|+|z_n-v_n|^\alpha)-{1\over n}\\
&\geq & \Dis -f_t(\omega)-C|y_n|-C|z_n|^\alpha-{1\over n}+(n-1)C(|y_n-u_n|+|z_n-v_n|^\alpha),
\end{array}
\end{equation}
which means that $\as$, in view of (ii),
$$
(n-1)C(|y_n-u_n|+|z_n-v_n|^\alpha)\leq
2(f_t(\omega)+C|y_n|+C|z_n|^\alpha)+{1\over n}
$$
and then
$$\limsup_{n\To \infty} nC(|y_n-u_n|+|z_n-v_n|^\alpha)<+\infty.\vspace{0.2cm}$$
Therefore, $\as$,
$$\Lim u_n=y_0,\ \ \Lim v_n=z_0.$$
Then, it follows from (21) and (H5) that, in view of Remark 6, $\as$,
$$\begin{array}{lll}
\Dis \liminf_{n\To \infty} g_n(\omega,t,y_n,z_n)&\geq &
\Dis \liminf_{n\To \infty} g(\omega,t,u_n,v_n)\\
&\geq & \Dis \liminf\limits_{(y,z)\To (y_0,z_0)}g(\omega,t,y,z)\\
&\geq & \Dis g(\omega,t,y_0,z_0).
\end{array}$$
On the other hand, from (20) and (18) we can also deduce that, $\as$,
$$\begin{array}{lll}
\Dis\limsup_{n\To \infty} g_n(\omega,t,y_n,z_n)&\leq &
\Dis \limsup_{n\To \infty} g(\omega,t,y_n,z_n)\\
&=& \Dis \lim\limits_{(y,z)\To (y_0^-, z_0)}g(\omega,t,y,z)\\
&=& \Dis g(\omega,t,y_0,z_0).
\end{array}$$
Hence, (iv) holds true, and the proof of Proposition 2 is complete.\hfill $\Box$\vspace{0.2cm}

{\bf Remark 8} Assume that the generator $g$ satisfies (H1b) and (H5). Similar argument to Proposition 2 yields that if we replace (20) with
$$
g^n(\omega,t,y,z)=\sup\limits_{(u,v)\in \Q^{1+d}}
(g(\omega,t,u,v)-nC(|y-u|+|z-v|^\alpha),
$$
then the conclusions of Proposition 2 hold also true for $g^n$, except that $g^n$ is non-increasing in $n$ and bigger than $g$, and that $y_0^-$ in (iv) is replaced with $y_0^+$.\vspace{0.2cm}

Now, we can begin the proof of Theorem 5.\vspace{0.1cm}

{\bf Proof of Theorem 5.}\ Suppose now that $\xi\in\LT$ and that (H1a) and (H5) hold for the generator $g$. For each $n\in\N$ and each $(\omega,t,y,z)\in \Omega\tim \T\tim \R\tim \R^d$, let $g_n(\omega,t,y,z)$ be defined in (20) and $h(\omega,t,y,z)$ be defined as follows
$$h(\omega,t,y,z):=f_t(\omega)+C(|y|+|z|^{\alpha}).$$
In view of (ii) and (iii) in Proposition 2, we know that both $h$ and $g_n$ are Lipschitz continuous in $y$ and $\alpha$-H\"{o}lder continuous in $z$, and $\as$, for each $n\geq 1$ and each $(y,z)\in \R\tim\R^d$, we have
$$g_n(\omega,t,y,z)\leq g_{n+1}(\omega,t,y,z)\leq g(\omega,t,y,z)\leq h(\omega,t,y,z).$$
Then, it follows from Theorem 1 in \citet{Fan10} that for each $n\geq 1$, both BSDE$(\xi,g_n)$ and BSDE$(\xi,h)$ have unique $L^1$ solutions, denoted, respectively,  by $(y^n_t,z^n_t)_{t\in \T}$ and $(\tilde{y}_t,\tilde{z}_t)_{t\in \T}$ for notational convenience. Furthermore, by Lemma 2 we also know that for each $t\in \T$,
\begin{equation}
  y^1_t\leq y^n_t\leq y^{n+1}_t\leq \tilde{y}_t\ \ \ \ps.
\end{equation}
We define $y_\cdot=\Lim y^n_\cdot$, then
\begin{equation}
\RE\ t\in[0,T],\ \ |y_t|\leq \sup\limits_{n\geq 1}|y^n_t|\leq |y^1_t|+|\tilde{y}_t|\ \ \ \ps.
\end{equation}

In the sequel, we will use the localization procedure again to construct the desired minimal solution. For each $k\geq 1$, introduce the following stopping time:
$$\tau_k=\inf\{t\in[0,T]: |y^1_t|+|\tilde{y}_t|+\int_0^t f_s \ {\rm d}s\geq k\}\wedge T.$$
Then $(y^n_{k}(t),z^n_{k}(t)):=(y^n_{t\wedge \tau_k },z^n_t 1_{t\leq \tau_k})$ solves the
following BSDE:
\begin{equation}
  y_{k}^n(t)=\xi^n_{k}+\int_t^T 1_{s\leq \tau_k}
  g_n(s,y_{k}^n(s),z_{k}^n(s))
  {\rm d}s-\int_t^T z_{k}^n(s)\cdot {\rm d}B_s,
\end{equation}
where $\xi^n_{k}=y^n_{\tau_k}$.

It is very important to observe that $y_{k}^n$ is nondecreasing in $n$ and that, from the definition of $\tau_k$ and inequality (23),
$$
\sup\limits_{n\geq 1}\sup\limits_{t\in[0,T]}
\|y_{k}^n(t)\|_\infty\leq k.
$$
Furthermore, by (ii) of Proposition 2 we have
$$|1_{s\leq \tau_k}g_n(s,y,z)|
\leq 1_{s\leq \tau_k}f_s+C+C|y|+C|z|,\ \ \RE\ n\geq 1.$$
Thus, in view of (iv) of Proposition 2 and the facts that $y^n_k$ is nondecreasing in $n$ and
$$\int_0^T \left(1_{s\leq \tau_k}f_s+C\right) {\rm d}s\leq k+CT,\vspace{0.1cm}$$
arguing as in the proof of Theorem 1, we can take the limit with respect to $n$ ($k$ being fixed) in (24) in the space $\s^2\times {\rm M}^2$. In particular, setting $y_k(t)=\sup_{n\geq 1}y^n_k(t)$, we know that $y_k(\cdot)$ is continuous and that there exists a process $z_{k}(t)\in {\M}^2$ such that $\Lim z^n_{k}(t)=z_{k}(t)$ in ${\M}^2$ and $(y_{k}(t),z_{k}(t))$ solves the BSDE
\begin{equation}
  y_{k}(t)=\xi_{k}+\int_t^T 1_{s\leq \tau_k}
  g(s,y_{k}(s)),z_{k}(s))
  {\rm d}s-\int_t^T z_{k}(s)\cdot {\rm d}B_s,
\end{equation}
where $\xi_{k}=\sup_{n\geq 1} y^n_{\tau_k}$.

Since $\tau_k\leq \tau_{k+1}$, it follows from the definitions of $y_{k}(\cdot),z_{k}(\cdot)$ and $y_\cdot$ that
$$y_{t\wedge\tau_k}=y_{k+1}(t\wedge\tau_k)=y_{k}(t)=\sup\limits_{n\geq 1}y^n_{t\wedge\tau_k},
\ \ z_{k+1}(t)1_{t\leq \tau_k}=z_{k}(t)=\lim\limits_{n\To \infty} z^n_t1_{t\leq \tau_k}.$$
Thus, since $y_{k}(\cdot)$ are continuous processes and moreover $\ps,\  \tau_k=T$ for $k$ large enough, we know that $y_\cdot$ is continuous on $[0,T]$. Then we define $z_\cdot$ on $(0,T)$ by setting
$$z_t=z_{k}(t),\ \ {\rm if}\ t\in(0,\tau_k),\vspace{0.3cm}$$
so that $z_t1_{t\leq \tau_k}=z_k(t)1_{t\leq \tau_k}=z_k(t)$ and (25) can be rewritten as
\begin{equation}
  y_{t\wedge\tau_k}=y_{\tau_k}
  +\int_{t\wedge\tau_k}^{\tau_k}
  g(s,y_s,z_s) {\rm d}s-\int_{t\wedge\tau_k}^{\tau_k}z_s\cdot {\rm d}B_s.
\end{equation}

Furthermore, we have
$$
\begin{array}{lll}
&& \Dis P \left(\int_0^T|z_s|^2{\rm d}s=\infty\right)\\
& =& \Dis P \left(\int_0^T|z_s|^2{\rm d}s=\infty,\tau_k=T\right)+ P \left(\int_0^T|z_s|^2{\rm
d}s=\infty,\tau_k<T\right)\\
&\leq & \Dis P \left(\int_0^{\tau_k}|z_{k}(s)|^2{\rm d}s=\infty\right)+
P\left(\tau_k<T\right),
\end{array}$$
and we deduce, since $\tau_k\uparrow T$, that
$$\int_0^T |z_s|^2{\rm d}s<\infty\ \ \ \ps.$$
Thus, note by (ii) and (iii) of Proposition 2 that for each $n\geq 1$, $g_n$ satisfies (H2) and (H4'') with $\mu=\gamma=nC$, and that it is smaller than $g$, letting $k\To \infty$ in (26) and arguing as in the proof of Theorem 1, we can deduce that $(y_t,z_t)$ is a minimal $L^1$ solution of BSDE$(\xi,g)$.

Finally, in view of Remark 8, using the same arguments as before we can prove the case of the maximal $L^1$ solution. Theorem 5 is then proved.\hfill $\Box$\vspace{0.2cm}

{\bf Remark 9}  Under the conditions (H1a) (resp. (H1b)) and (H5), it is uncertain whether the $L^1$ solution of BSDE$(\xi,g)$ is unique or not, an counterexample can be found in \citet{Jia08}.\vspace{0.2cm}

With Theorem 5 in hand, using the same arguments as in Theorems 2-4 and Remarks 4-5 and noticing the fact that (H1') can imply not only (H1a) but also (H1b), we can obtain the following Theorems 6-8.\vspace{0.1cm}

{\bf Theorem 6} (Comparison theorem on the minimal (resp. maximal) $L^1$ solution)\ Assume that $\xi,\xi'\in \LT$ and that both $g$ and $g'$ satisfy (H1a) (resp. (H1b)) and (H5). Let $(y^{g}_t(\xi),z^{g}_t(\xi))_{t\in \T}$ and $(y^{g'}_t(\xi'),z^{g'}_t(\xi'))_{t\in \T}$ be, respectively, the minimal (resp. maximal) $L^1$ solution to BSDE$(\xi,g)$ and BSDE$(\xi',g')$ by Theorem 5. If $\ps,\ \xi\leq \xi'$ and $\as$, for each $(y,z)\in \R\tim \R^d$, $g(t,y,z)\leq g'(t,y,z)$, then for each $t\in\T$, we have
$$y^{g}_t(\xi)\leq y^{g'}_t(\xi')\ \ \ \ps.\vspace{-0.1cm}$$

{\bf Theorem 7} \ (Levi type theorem on the minimal (resp. maximal) $L^1$ solution)\ Assume that $\xi_n,\xi\in \LT$ for each $n\geq 1$ and that $g$ satisfies (H1a) (resp. (H1b)) and (H5). Let $(y^{g}_t(\xi_n),z^{g}_t(\xi_n))_{t\in \T}$ and $(y^{g}_t(\xi),z^{g}_t(\xi))_{t\in \T}$ be, respectively, the minimal (resp. maximal) $L^1$ solution of BSDE$(\xi_n,g)$ and BSDE$(\xi,g)$ by Theorem 5. If $\ps$, $\xi_n\uparrow\xi$ (resp. $\xi_n\downarrow\xi$), then for each $t\in\T$,
$$\Lim y^{g}_t(\xi_n)=y^{g}_t(\xi)\ \ \ \ps.\vspace{-0.2cm}$$

{\bf Theorem 8}\ (Lebesgue type theorem on the $L^1$ solution)\ Assume that $\xi_n,\xi\in \LT$ for each $n\geq 1$ and that $g$ satisfies (H1') and (H5). Assume further that BSDE$(\xi,g)$ has a unique $L^1$ solution $(y^{g}_t(\xi),z^{g}_t(\xi))_{t\in \T}$. Let $(y^{g}_t(\xi_n),z^{g}_t(\xi_n))_{t\in \T}$ be any of $L^1$ solutions of BSDE$(\xi_n,g)$ by Theorem 5. If $\ps,\ \xi_n\To\xi$ as $n\To\infty$ and $\ps,\ |\xi_n|\leq \eta$ with $\E[|\eta|]<+\infty$, then for each $t\in\T$,
$$\Lim y^{g}_t(\xi_n)=y^{g}_t(\xi)\ \ \ \ps.$$

\section{A general comparison theorem on $L^1$ solutions }

In this section, under the assumptions that $g$ is weakly monotonic in $y$ and uniformly continuous in $z$ as well as it has a stronger sublinear growth in $z$, we will establish a general comparison theorem on $L^1$ solutions of the BSDEs. Let us introduce the following assumptions taken from \citet{Fan12a}\vspace{0.2cm}:

(H2') $g$ is weakly monotonic in $y$, i.e., there exists a
nondecreasing concave function $\rho(\cdot)$ from ${\R}_+$ to itself with $\rho(0)=0$, $\rho(u)>0$ for $u>0$ and $\int_{0^+}{1\over \rho(u)}\ {\rm d}u=+\infty$ such that $\as,$
$$
(g(\omega,t,y_1,z)-g(\omega,t,y_2,z))\ {\rm sgn}(y_1-y_2) \leq\rho(|y_1-y_2|),\ \ \RE\ y_1,y_2,z;\ \ \
$$

(H4*) $g$ is uniformly continuous in $z$ uniformly with respect to $(\omega,t,y)$, i.e., there exists a
continuous, nondecreasing function $\phi(\cdot)$ from ${\R}_+$ to itself with linear growth and satisfying $\phi (0)=0$ such that $\as,$
$$\RE\ y,z_1,z_2,\ \
|g(\omega,t,y,z_1)-g(\omega,t,y,z_2)|\leq\phi(|z_1-z_2|).
$$

{\bf Remark 10} It is clear that (H2') and (H4*) are, respectively, weaker than (H2) and (H4'').\vspace{0.2cm}

Using the similar arguments to Theorem 1 in \citet{Fan12a} together with the stopping time technique, we can obtain the following Proposition 3. It is a slight generalization of Theorem 1 in \citet{Fan12a}, where only is the $L^2$ solution to BSDEs investigated.\vspace{0.1cm}

{\bf Proposition 3} (Comparison theorem)\ Let $g$ and $g'$ be two generators of BSDEs, and let $(y_\cdot,z_\cdot)$ and $(y'_\cdot,z'_\cdot)$ be, respectively, a solution to BSDE$(\xi,g)$ and BSDE$(\xi',g')$. Assume that $\ps,\ \xi\leq \xi'$, $g$ satisfies (H2') and (H4*), and $\as,\ g(t,y'_t,z'_t)\leq g'(t,y'_t,z'_t)$ (or $g'$ satisfies (H2') and (H4*), and $\as,\ g(t,y_t,z_t)\leq g'(t,y_t,z_t)$). If $(y_\cdot-y'_\cdot)^+$ belongs to $\s$, then for each $t\in\T$, we have
$$y_t\leq y'_t\ \ \ \ps.\vspace{0.1cm}$$

By virtue of the above Proposition 3, we can prove the following comparison theorem on the $L^1$ solutions of BSDEs, which improves Proposition 1 in \citet{Fan10} and Proposition 2 in \citet{Xiao12}.\vspace{0.1cm}

{\bf Theorem 9} (Comparison theorem on the $L^1$ solution)\ Let $g$ and $g'$ be two generators of BSDEs, and let $(y_\cdot,z_\cdot)$ and $(y'_\cdot,z'_\cdot)$ be, respectively, a $L^1$ solution to BSDE$(\xi,g)$ and BSDE$(\xi',g')$. If $\ps,\ \xi\leq \xi'$, $g$ satisfies (H2'), (H4') and (H4*), and $\as,\ g(t,y'_t,z'_t)\leq g'(t,y'_t,z'_t)$ (or $g'$ satisfies (H2'), (H4') and (H4*), and $\as,\ g(t,y_t,z_t)\leq g'(t,y_t,z_t)$), then for each $t\in\T$,
$$y_t\leq y'_t\ \ \ \ps.$$

{\bf Proof.}\ It follows from Proposition 3 that we need only to show that $(y_\cdot-y'_\cdot)^+$ belongs to $\s$ under the assumptions of Theorem 9.

Now, we assume that $\ps,\ \xi\leq \xi'$, $g$ satisfies (H2'), (H4') and (H4*), and $\as,\ g(t,y'_t,z'_t)\leq g'(t,y'_t,z'_t)$. The same arguments as follows can prove the another case. Let us fix $k\in \N$ and
denote the stopping time
$$\tau_k:=\inf \left\{t\in\T: \int_0^t\left(|z_s|^2+|z'_s|^2
\right)\ {\rm d}s\geq k\right\}\wedge T.
$$
Tanaka's formula leads to the equation, setting $\hat{y}_t=y_t-y'_t,\
\hat{z}_t=z_t-z'_t$,
$$
\Dis \hat{y}_{t\wedge \tau_k}^+\leq \hat{y}_{\tau_k}^+
+\int_{t\wedge \tau_k}^{\tau_k}
 1_{\hat{y}_s> 0}[g(s,y_s,z_s)-g'(s,y'_s,z'_s)]
 \ {\rm d}s-\int_{t\wedge \tau_k}^{\tau_k}1_{\hat{y}_s> 0}
\hat{z}_s\cdot {\rm d}B_s.$$
Since $g(s,y'_s,z'_s)-g'(s,y'_s,z'_s)$ is
non-positive, we have
$$
\begin{array}{lll}
g(s,y_s,z_s)-g'(s,y'_s,z'_s)&=& g(s,y_s,z_s)-g(s,y'_s,z'_s)+g(s,y'_s,z'_s)-g'(s,y'_s,z'_s)\\
&\leq &g(s,y_s,z_s)-g(s,y'_s,z_s)+g(s,y'_s,z_s)-g(s,y'_s,z'_s),
\end{array}
$$
and we deduce, using the assumptions (H2') and (H4') of $g$, that
$$
1_{\hat{y}_s> 0}[g(s,y_s,z_s)-g'(s,y'_s,z'_s)]\leq
\rho(\hat{y}_s^+)+2\lambda(f_s+|y'_s|+|z_s|+|z'_s|)^{\alpha}.
$$
Thus, we get that
$$
\Dis \hat{y}_{t\wedge \tau_k}^ +\leq \hat{y}_{\tau_k}^+
+\int_{t\wedge \tau_k}^ { \tau_k}\left[\rho(\hat{y}_s^+)+2\lambda(f_s+|y'_s|+
|z_s|+|z'_s|)^{\alpha}\right]{\rm d}s-\int_{t\wedge
\tau_k}^ { \tau_k}1_{\hat{y}_s> 0}\hat{z}_s\cdot {\rm d}B_s,
$$
and then that
\begin{equation}
\Dis\hat{y}_{t\wedge \tau_k}^ + \leq
\Dis\E\left[\hat{y}_{\tau_k}^+ +\left.\int_{t\wedge
\tau_k}^{\tau_k}
\left[\rho(\hat{y}_s^+)+2\lambda(f_s+|y'_s|+|z_s|+
|z'_s|)^{\alpha}\right] {\rm d}s\right|\F_{t}\right].
\end{equation}
Furthermore, since $\rho(\cdot)$ is a nondecreasing concave function and $\rho(0)=0$, we can find a pair of positive constants $a$ and $b$ such that
\begin{equation}
\rho(u)\leq a+bu,\ \ \RE\  u\geq 0.
\end{equation}
Then, since both $(y_\cdot,z_\cdot)$ and $(y'_\cdot,z'_\cdot)$ are $L^1$ solutions, we can send $k$ to $\infty$ in (27) and use the Lebesgue dominated convergence theorem, in view of $\xi\leq \xi'$,  $\tau_k\To T$ as $k\To \infty$ and (28), to get that, for each $t\in \T$,
$$
\begin{array}{lll}
\Dis \hat{y}_t^+ &\leq & \Dis 2\lambda\E\left[\left.\int_0^T
(f_s+|y'_s|+|z_s|+|z'_s|)^{\alpha}{\rm
d}s\right|\F_t\right]+\E\left[\left.\int_t^T \rho(\hat{y}_s^+){\rm
d}s\right|\F_t\right]\\
&\leq &\Dis aT+2\lambda \E\left[\left.\int_0^T
(f_s+|y'_s|+|z_s|+|z'_s|)^{\alpha}{\rm d}s\right|\F_t\right]+b\int_t^T
\E\left[\left.\hat{y}_s^+\right|\F_t\right]{\rm d}s,
\end{array}\vspace{0.2cm}$$
and then for each $r\in [t,T]$,
$$\E\left[\left.\hat{y}_r^+\right|\F_t\right]\leq aT+2\lambda
\E\left[\left.\int_0^T (f_s+|y'_s|+|z_s|+
|z'_s|)^{\alpha}{\rm d}s\right|\F_t\right]+b\int_r^T
\E\left[\left.\hat{y}_s^+\right|\F_t\right]{\rm d}s.\vspace{0.2cm}$$
Gronwall's inequality yields that for each $r\in [t,T]$,
$$\E\left[\left.\hat{y}_r^+\right|\F_t\right]\leq
\left(aT+2\lambda\E\left[\left.\int_0^T (f_s+|y'_s|+|z_s|+|z'_s|)^{\alpha}{\rm
d}s\right|\F_t\right]\right)\cdot e^{b(T-r)},$$ form which, by letting $r=t$, we have
$$\hat{y}_t^+\leq \left(aT+2\lambda\E\left[\left.\int_0^T
(f_s+|y'_s|+|z_s|+|z'_s|)^{\alpha}{\rm d}s\right|\F_t\right]\right)\cdot e^{bT}.$$
Finally, taking supremum over $t$ and then taking expectation in both sides of the above inequality follows that, by virtue of Doob's inequality, H$\ddot{{\rm o}}$lder's inequality and the fact that both $(y_\cdot,z_\cdot)$ and $(y'_\cdot,z'_\cdot)$ are $L^1$ solutions,
$$
\begin{array}{lll}
\Dis \E[\sup_{t\in\T}|\hat{y}_t^+|^{\beta/\alpha}]&\leq & \Dis K\left(1+\E\left[\left(\int_0^T f_s {\rm d}s
\right)^{\beta}\right]+\E[\sup_{t\in\T}|y'_t|^{\beta}]\right.\\
&& \ \ \Dis \left.+\E\left[\left(\int_0^T |z_s|^2{\rm d}s
\right)^{\beta/2}\right]+
\E\left[\left(\int_0^T |z'_s|^2{\rm d}s
\right)^{\beta/2}\right]\right)\\
&<&+\infty,
\end{array}
$$
where $\beta$ is any constant which belongs to $(\alpha,1)$, and $K$ is a constant depending only on $(a,b,T,\lambda,\alpha,\beta)$. That is to say,
$\hat y_\cdot^+=(y_\cdot-y'_\cdot)^+ \in {\s}$.
Then the proof of Theorem 9 is completed. \hfill $\Box$\vspace{0.2cm}

Combining Theorem 9 with Theorem 1, in view of Remarks 1 and 10, we can obtain the following existence and uniqueness result.\vspace{0.1cm}

{\bf Theorem 10} (Existence and uniqueness theorem on the $L^1$ solution)\ Assume that the generator $g$ satisfies (H1)-(H3), (H4') and (H4*). Then for each $\xi\in\LT$, BSDE$(\xi,g)$ has a unique $L^1$ solution.\vspace{0.2cm}

{\bf Remark 11}\ Compared with the one-dimensional versions of Theorems 6.2 and 6.3 in \citet{Bri03}, we can see that the Lipschitz continuity assumption of $g$ in $z$ employed in \citet{Bri03} is weakened to the uniform continuity assumption (H4*) here. \vspace{0.2cm}

{\bf Example 3}\ For each $(\omega,t,y,z)\in \Omega\times\T\times\R\times\R^d$, let
$$g(\omega,t,y,z)=|B_t(\omega)|^2\cdot{\text e}^{-y}+\sqrt{1+|y|+|z|}+\sqrt[3]{|z|}+{1\over \sqrt{|t-T/2|}}1_{t\neq T/2}.\vspace{-0.1cm}$$
It is not hard to check that this $g$ satisfies assumptions (H1)-(H3), (H4') and (H4*) with $\mu=1,\lambda=2,\alpha=1/2$ and $f_t(\omega)\equiv 1$. It then follows from Theorem 10 that for each $\xi\in\LT$, BSDE$(\xi,g)$ has a unique $L^1$ solution.

It should be especially pointed out that this generator $g$ has a general growth in the variable $y$, it is uniformly continuous with respect to the variable $z$, but it is neither Lipschitz continuous nor H\"{o}lder continuous in $z$. Then, the existence and uniqueness result of $L^1$ solutions to BSDE$(\xi,g)$ with $\xi\in\LT$ can not be obtained by any existing results.
\vspace{0.4cm}




\end{document}